\documentclass[11pt]{article}
\usepackage{amsmath}
\usepackage{amssymb}
\usepackage{amsfonts}
\usepackage{mathrsfs}
\usepackage{shortvrb,psfrag}
\usepackage{enumerate}
\usepackage{color}
\usepackage[dvips]{graphicx}
\usepackage{pdfsync}

\newtheorem{defn}{Definition}
[section]
\newtheorem{lemma}{Lemma}
[section]

\newtheorem{thm}{Theorem}
[section]
\newtheorem{cor}{Corollary}
[section]
\newtheorem{rem}{Remark}
[section]

\newcommand{\la}{\lambda}

\newcommand{\ptl}{\partial}

\newcommand{\ra}{\rightarrow}
\newcommand{\lra}{\longrightarrow}

\def\a{\alpha}

\def\R{\Bbb R}

\def\ba{\begin{array}}
\def\ea{\end{array}}

\def\dag{\dagger}
\def\part{\partial}
\def\fr{\frac}

 \newcommand{\beq}{\begin{equation}}
\newcommand{\eeq}{\end{equation}}
\newcommand{\ben}{\begin{eqnarray}}
\newcommand{\een}{\end{eqnarray}}
\newcommand{\beno}{\begin{eqnarray*}}
\newcommand{\eeno}{\end{eqnarray*}}

\newcommand{\ls}{\setlength{\baselineskip}{18pt}
                      \setlength{\parskip}{3mm} }

\title {$C^{\a}$ regularity of weak solutions of  non-homogenous ultraparabolic equations with drift terms}
\author{Wendong Wang$^\dag$\, Liqun Zhang$^\ddag$\\[2mm]
{\small $^\dag$School of  Mathematical Sciences, Dalian University of Technology, Dalian 116024, P.R. China}\\[1mm]
{\small E-mail: wendong@dlut.edu.cn}\\[1mm]
{\small $^\ddag$Institute of Mathematics, AMSS, School of Mathematical Sciences, UCAS,}\\[1mm]
{\small Chinese Academy of Sciences,
Beijing 100190, P.R. China}\\[1mm]
{\small E-mail: lqzhang@math.ac.cn}}

\begin{document}

\maketitle

\begin{abstract}
Consider a class of non-homogenous ultraparabolic differential equations
with drift terms or lower order terms arising from some physical models, and we prove that weak solutions are H\"{o}lder continuous, which also generalizes the classic results of  parabolic equations of second order. The main ingredients are a type of weak Poincar\'{e} inequality satisfied by non-negative weak sub-solutions and Moser iteration.
\end{abstract}

{\small {\bf Keywords:} Ultraparabolic equations, Moser iteration, Poincar\'{e} inequality, $C^{\a}$ regularity }

% Smoothness and regularity of solutions, Hypoelliptic equations

{\bf 1991 Mathematics Subject Classification.} 35K70, 35H10, 35B65

\pagenumbering{arabic}
%%%%%%%%%%%%%%%%%%%%%%%%%%%%%%%%%%%%%%%%%%%%%%%%%%%%%%%%%%%%%%%%%%%%%%%%%%%%%%%%
\setcounter{equation}{0}
\section{Introduction} \ls \noindent

Ultraparabolic equations form a class of degenerated parabolic equations, which come from kinetic equations, diffusion process, Asian options and so on.
One of the typical examples of ultraparabolic equations is the
following equation of Kolmogorov type:
\ben\label{eq:KOL}
\displaystyle \fr{\partial \,u}{\partial \,t}+y \fr{\partial
\,u}{\partial \,x}-\fr{\partial^2 \,u}{\partial
\,y^2}=0,
\een
which is a strongly degenerated parabolic equation.

On the other hand, the results obtained by Pascucci-Polidoro \cite{PP} proved that the
Moser iteration still works for a class of ultraparabolic
equations with measurable coefficients; see also, Cinti-Pascucci-Polidoro \cite{CPP} dealing with the non-homogenous case and Cinti-Polidoro \cite{CP} for  general hypoelliptic equations. Their results show that for
a non-negative sub-solution $u$ of (\ref{eq:KOL}), the $L^{\infty}$ norm of
$u$ is bounded by the $L^p$ norm ($p \ge 1$). This is a very
important step to the final regularity of solutions of the
ultraparabolic equations. Based on these estimates for weak solutions, Zhang \cite{Zhang} for the homogenous ultraparabolic
equations, Wang-Zhang \cite{Wang-Zhang2009} for the non-homogenous case and Wang-Zhang \cite{Wang-Zhang2011} for general hypoelliptic equations can prove the H\"{o}lder regularity of weak solutions with the help of De Giorgi-Nash-Moser iteration by exploring a weak Poincar\'{e} inequality; see also Wang-Zhang \cite{Wang-Zhang3} for other different ultraparabolic parabolic cases. Note that the above progress is based on the weak solutions
without lower order terms(for example $c(x)u$ or $f$). However, when lower order terms exist, their integrability in different critical Sobolev spaces might change the relevant weak Sobolev and Poincar\'{e} inequalities, which cause some difficulties in  De Giorgi-Nash-Moser iterations, since one can not use the Log form auxiliary function directly (see, for example,  Lemma 3.2 from below).

In this paper, we are concerned with the $C^{\a}$ regularity of
solutions of more general ultraparabolic equations and consider the following non-homogeneous Kolmogorov-Fokker-Planck type
operator on $\mathbb{R}^{N+1}$:
\ben\label{eq:general equ}
Lu =\sum_{i=1}^{m_0}b'_i(x,t)\partial_iu+c(x,t)u+f(x,t),
\een
where $L$ is given by:
\beno
Lu \equiv
\sum_{i,j=1}^{m_0}{\ptl_{i}(a_{ij}(x,t)\ptl_{j}\,u
)}+\sum_{i,j=1}^N b_{ij}x_i {\ptl_{j}\, u }- {\ptl_t
\,u}
\eeno
Here $(x,t)\in \mathbb{R}^{N+1}$, $1\leq m_0\leq N$, $\ptl_{x_j}=\ptl_{j}$ and $b_{ij}$ is
constant for every $i,j= 1,\cdots, N$. Let $A=(a_{ij})_{N\times N},$
where $a_{ij}=0,$ if $i>m_0$ or $j>m_0.$ Moreover, $b'(x,t)\in \mathbb{R}^{m_0},$ $c(x,t),f(x,t)\in \mathbb{R}$ are measurable functions. We make the following
assumptions on the coefficients of $L$:

$({\bf H_1})$  $a_{ij}=a_{ji} \in L^{\infty} (\mathbb{R}^{N+1})$ and there
exists a $\lambda >0$ such that
$$
\fr{1}{\lambda}\sum_{i=1}^{m_0}\xi_i^2 \leq \sum_{i,j=1}^{m_0}
a_{ij}(x,t)\xi_i \xi_j \leq {\lambda}\sum_{i=1}^{m_0}\xi_i^2
$$
for every $(x,t)\in \mathbb{R}^{N+1}$, and $\xi \in \mathbb{R}^{m_0}$.

$({\bf H_2})$ The matrix $B=(b_{ij})_{N \times N}$ has the form
$$\left(
\begin{array}{ccccc}
\ast & {B_1} & 0 & \cdots & 0  \\
\ast &  \ast  & {B_2} & \cdots & 0 \\
\vdots & \vdots & \vdots & \ddots &\vdots \\
\ast &  \ast & \ast & \cdots & {B_d} \\
\ast &  \ast & \ast & \cdots & \ast
\end{array}
\right)
$$
where $B_k$ is a matrix $m_{k-1}\times m_{k}$ with rank $m_k$,
$m_0\geq m_1\geq \cdots\geq m_d$ with $m_0+m_1+\cdots+m_d=N$ and $\ast$ denotes any real number.

$({\bf H_3})$: Let $c,f,b'(x,t)$ satisfy the conditions:
\ben\label{eq:b' condi}
c,f\in L^q(\Omega)~{\rm with}~q>\frac{Q+2}{2},\quad b'\in L^{Q+2}(\Omega).
%b'\in L_t^sL_x^p(\Omega),\quad \frac{Q}{p}+\frac{2}{s}=1,\quad 1\leq s <\infty.
\een

The requirements of matrix $B$ in $(H_2)$ ensure that the operator
$L$ with the constant $a_{ij}$ satisfies
H\"ormander's hypoellipticity condition. We refer to \cite{CPP} for more details on non-homogeneous
Kolmogorov-Fokker-Planck type operator on $\mathbb{R}^{N+1}$.

Schauder estimates for the solutions to (\ref{eq:general equ}) have been obtained for example,
in \cite{Lun,Man97,FP}. In addition, the H\"{o}lder regularity of weak solutions
have been proved by Bramanti-Cerutti-Manfredini \cite{BCM} and Polidoro-Ragusa \cite{PR} by assuming a weak continuity on the coefficient
$a_{ij}$; see also Angiuli-Lorenzi \cite{AL} for higher derivatives estimates by assuming the coefficients are H\"{o}lder continuous or smooth. It is quite interesting to establish whether the weak solution has
H\"older regularity under the assumption $(H_1)$ on $a_{ij}$. One of
the approach to the H\"older estimates is to obtain the Harnack type
inequality. In the case of elliptic equations with measurable
coefficients, the Harnack inequality was obtained by Moser \cite{Moser61} via
an estimate of BMO functions due to the John-Nirenberg inequality
together with the Moser iteration. Moser \cite{Moser64} also
obtained the Harnack inequality for parabolic equations with
measurable coefficients by generalizing the John-Nirenberg estimates
to the parabolic case. Another approach to the H\"older regularity is
given by Kruzhkov \cite{Kru63,Kru64} based on the Moser iteration to
obtain a local priori estimates, which provides a short proof for
the parabolic equations. Earlier, De Giorgi \cite{De} developed an approach to obtain
the H\"older regularity for elliptic equations. Nash \cite{Nash} also introduced another technique
relying on the Nash inequality and obtained the
H\"older regularity.

Our main idea is to establish a type of weak ${\rm Poincar\acute{e}}$ type inequality for non-negative
weak sub-solutions used by the authors in \cite{Zhang}, \cite{Wang-Zhang2009} and \cite{Wang-Zhang2011} in order to extend the main
results in \cite{Zhang}, \cite{Wang-Zhang2009} and \cite{Wang-Zhang2011} to operators with lower order terms. Then by using Kruzhkov's method of level sets we can obtain local a
priori estimates which implies the H\"older estimates for
ultraparabolic equation (\ref{eq:general equ}).

Next, we give a detailed definition of weak solution. Let $D_{m_0}$ be the gradient with respect to the variables $x_1,
x_2,\cdots, x_{m_0}$. As in \cite{CPP}, we write
$$Y=\sum_{i,j=1}^N b_{ij}x_i {\ptl_{x_j}}- {\ptl_t}.$$
\begin{defn}
We say that $u$ is a {\it weak solution} of (\ref{eq:general equ}) in a domain $\Omega\subset \mathbb{R}^{N+1}$ if it satisfies (\ref{eq:general equ}) in the
$H^1$ sense, that is for any $\phi \in C^{\infty}_0(\Omega)$, there holds
\ben\label{eq:weak solution}
\int_{\Omega} \phi Yu-(Du)^T AD\phi =\int_{\Omega} \phi(b'\cdot D_{m_0}u+ cu+f),
\een
and $u$, $D_{m_0}u$, $Yu, b', c, f \in L^2_{\rm loc}(\Omega)$.
\end{defn}

Similarly, we can define the {\it weak sub-solutions (super-solutions)} of (\ref{eq:general equ}) in a domain $\Omega\subset \mathbb{R}^{N+1}$, if $u$, $D_{m_0}u$, $Yu, b', c, f \in L^2_{\rm loc}(\Omega)$, and  for any nonnegative $\phi \in C^{\infty}_0(\Omega)$, there holds
\ben\label{eq:weak subsolution}
\int_{\Omega} \phi Yu-(Du)^T AD\phi \geq(\leq)\int_{\Omega} \phi(b'\cdot D_{m_0}u+ cu+f).
\een

One of the important feature of equation (\ref{eq:general equ}) is that the
fundamental solution can be written explicitly if the coefficients
$a_{ij}$, $b_j'$ and $c$ are constant (cf. \cite{Kup,LanP}). Besides, there are some
geometric and algebraic structures in the space $\mathbb{R}^{N+1}$ induced by
the constant matrix $B$.

Recently in \cite{GIMV},  the authors consider the Laudau equation,
\beno
\partial_t f +v\cdot \nabla_x f=\nabla_v\cdot(A[f]\nabla f+B[f]f)
\eeno
and obtain H\"{o}lder continuity and the Harnack inequality of weak solution by De Giorgi's method with bounded $A[f], B[f]$ by using an estimate
of \cite{Bou} for regularity gain in $x$ direction. It seems to be difficult to obtain the similar Harnack inequality in the general
 case (\ref{eq:general equ}), since Theorem 1.3 in \cite{Bou} can't be applied directly when the new first order terms through
  multiple Lie operations arise (for example, see the general form of $B$ in $({\bf H_2})$), which are not considered in \cite{Bou}.
 %Besides, we refer to \cite{IS} for regularity of weak solution of the Boltzman equation.

Now let us introduce some basic properties of hypoelliptic operator.
Let $E(\tau)=\rm{exp}(-\tau B^T)$. For $(x,t), (\xi,\tau) \in
\mathbb{R}^{N+1}$, set
$$(x,t)\circ (\xi,\tau)=(\xi+E(\tau)x,t+\tau),$$ then $(\mathbb{R}^{N+1}, \circ)$
is a Lie group with identity element $(0,0)$, and the inverse of an
element is $(x,t)^{-1}=(-E(-t)x,-t)$. The left translation by
$(\xi,\tau)$ given by
$$(x,t)\mapsto (\xi,\tau)\circ (x,t),$$
is a invariant translation to operator $L$ when coefficient $a_{ij}$
is constant.
The dilation associated to the operator $L$ with constant
coefficient $a_{ij}$ is given by
$$
\delta_t=diag (t I_{m_0},t^3 I_{m_1},\cdots,t^{2d+1}I_{m_d},t^2),
$$
where $I_{m_k}$ denotes the $m_k\times m_k$ identity matrix and $t$ is
a positive parameter. We define
$$
\tilde{D}_t=diag(t I_{m_0},t^3 I_{m_1},\cdots,t^{2d+1}I_{m_d}),
$$
and denote
$$Q=m_0+3m_1+\cdots+(2d+1)m_d,$$
then the number $Q+2$ is usually called the homogeneous dimension of
$(\mathbb{R}^{N+1},\circ)$ with respect to the dilation $\delta_t$. Note that $Q+2$ is a natural scaling for the group $(\mathbb{R}^{N+1}, \circ)$. Specially, $Q$ can be seen as the spatial dimension and the coefficients $c,f,b'(x,t)$ need to satisfy the similar conditions as in the Euclidean space (see $({\bf H_3})$).

Our main result is the following theorem.
\begin{thm}\label{thm:main}
Under the assumptions $({\bf H_1}-{\bf H_3})$, weak solutions of
(\ref{eq:general equ}) are H\"older continuous in the interior of $\Omega$.
\end{thm}

\begin{rem}
i)When $m_0=N$, we have $Q=N$, and at this moment (\ref{eq:general equ}) becomes a parabolic equation with lower order terms. The requirements of coefficients in
$({\bf H_1}-{\bf H_3})$ agree with the sharp form of the inequality for parabolic equations.\\
ii) The above assumptions on $b',c,f$
are due to the following embedding inequality. The conditions of $D_{m_0}u\in L^2_{loc}(\Omega) $, $u\in L_t^{\infty}L^2_x$ locally derived from (\ref{eq:weak solution}) via the basic energy estimate, and the non-negative weak sub-solution $u$ imply that
$u\in L^{\frac{2Q+4}{Q}}_{loc}(\Omega)$ (see (\ref{eq:Sobolev estimate-w}) from the following Lemma \ref{lem:Sobolev estimate} ).
%iii) The condition $\|f\|_{L^{Q+2}({\cal B}^-_{2}(x_0,t_0))}<\infty$ seems to be stronger. However, this is related to the weak Sobolev inequality. When $f\in L^q$ with $q>\frac{Q+2}{2}$, we only obtained the following embedding inequality
%(see Lemma \ref{lem:Sobolev estimate})
%for $w=u^p$ with the positive integer $p> 1$,
%\beno
%\|\phi w\|_{L^{2k}({\cal B}^-_{\rho}(x_0,t_0))}&\leq& \frac{C}{r-\rho}(p^{\frac{1}{1-\beta}}\|w\|_{L^2({\cal B}^-_{r}(x_0,t_0))}+ \|D_{m_0}w\|_{L^2({\cal B}^-_{r}(x_0,t_0))})\nonumber\\
%&&+ r^{\max\{p(Q+2)(\frac{1}{Q+2}-\frac1q+\frac{1}{2p}),0\}}\|f\|_{L^q({\cal B}^-_{r}(x_0,t_0))}^p
%\eeno
%where $\beta=\frac{Q+2}{q}-1\in (0,1)$ with $q>(Q+2)/2$. The power of $r$ in the front of $f$
%is not enough for the $L^{\infty}$ iterative scheme.
\end{rem}

%Applying the above result to the equations (\ref{eq:Prandtl}), since there exists a global weak solution of (\ref{eq:P}) obtained in \cite{XZ}, we have the following conclusion.
%
%\begin{cor}\label{cor:main}{\rm(Interior regularity criterion for 2D Prandtl equation)}
%Under the monotone class assumptions (\ref{eq:positivity})-(\ref{eq:monotone}) and the favorable pressure condition (\ref{eq:favourable}), the global weak solution of (\ref{eq:P}) obtained in \cite{XZ} is
%locally H\"{o}lder continuous if the additional conditions $\partial_tu, \partial_tU, \partial_xU\in L^6_{loc}$ are put.
%\end{cor}
%
%
%Recently, Li-Wu-Xu in \cite{LWX} obtained Gevrey class smoothing effect for the Prandtl equation with monotonic data by assuming $e^{cy}\partial_xu \in L^{\infty}(0,T;L^7(\mathbb{R}^2_+))$ and $e^{cy}\partial_{xy}u \in L^{2}(0,T;L^7(\mathbb{R}^2_+))$.

The paper is organized as follows. In section 2, we introduce the   Lie group structure on the ultraparabolic operator and properties of fundamental solution. Section 3 is devoted to obtaining some technical lemmas in proof of Theorem \ref{thm:main}, including level set estimate with the G-function method, weak Sobolev inequality, and weak Poincar\'{e} inequality. In section 4, we complete the proof of  Theorem \ref{thm:main}. The last section is an introduction to the G-function.

\setcounter{equation}{0}
\section{Preliminary Results on Lie Groups}

Given $r>0$,
the norm in $(\mathbb{R}^{N+1},\circ)$, related to the group of translations and
dilation to the equation is defined by $$||(x,t)||=r,$$ if $r$ is
the unique positive solution to the equation
$$
\fr{x_1^2}{r^{2\a_1}}+\fr{x_2^2}{r^{2\a_2}}+\cdots+\fr{x_N^2}{r^{2\a_N}}
+\fr{t^2}{r^4}=1,
$$
where $(x,t) \in \mathbb{R}^{N+1}\setminus \{0\}$ and
$$
\a_1=\cdots=\a_{m_0}=1, \quad
\a_{m_0+1}=\cdots=\a_{m_0+m_1}=3,\cdots,
$$
$$
\a_{m_0+\cdots+m_{d-1}+1}=\cdots=\a_N=2d+1.
$$
And $||0||=0$. The balls at a point $(x_0,t_0)$ is defined by
$${\cal B}_r(x_0,t_0)=\{(x,t)|\quad ||(x_0,t_0)^{-1}\circ (x,t)||\leq r\},$$
and
$${\cal B}^-_r(x_0,t_0)={\cal B}_r(x_0,t_0)\cap\{t<t_0\}.$$
For convenience, we sometimes use the cube replace the balls. The
cube at point $0$ is given by
$$
{\cal C}_r(0)=\{(x,t)|\quad |t|\leq r^2,\quad |x_1|\leq r^{\a_1},
\cdots, |x_N|\leq r^{\a_N}\}.
$$
and ${\cal C}_r(0)$
It is easy to see that there exists a constant $\Lambda$ such that
$$
{\cal C}_{\fr r \Lambda}(0)\subset{\cal B}_r(0)\subset{\cal
C}_{\Lambda r}(0),
$$
where $\Lambda$ only depends on $N$.

Let $A_0$ be the matrix $(a_{ij})_{N \times N}$. When $A_0$ is constant matrix and has the form
$$
A_0= \left(
\begin{array}{cc}
I_{m_0\times m_0}& 0  \\
0 & 0
\end{array}
\right)
$$
then let
$$
{\mathcal C}(t)\equiv\int_0^t E(s)A_0 E^T(s)ds,
$$
which is positive when $t>0.$ Moreover, we let the operator $L_1$ take the form
$$
L_1={\rm div} (A_0D)+Y,
$$
whose fundamental solution $\Gamma_1(\cdot,\zeta)$ with a pole at
$\zeta=0\in \mathbb{R}^{N+1}$ is:
$$
\Gamma_1(z,\zeta)=\Gamma_1(\zeta^{-1}\circ z,0), \qquad z, \zeta \in
\mathbb{R}^{N+1},\quad z \neq \zeta,
$$
where $z=(x,t)$. Especially, $\Gamma_1(z,0)$ can be written down explicitly
\ben\label{eq:Gamma1}
\Gamma_1(z,0)= \left\{
\begin{array}{lll}{\frac{(4\,\pi)^{-\frac{N}{2}}}{\sqrt{\det {\mathcal C}(t)}}\exp(-\frac{1}{4}\langle {\mathcal C}
^{-1}(t)x,x\rangle -t \, {\rm tr}(B))} & {\rm if} \quad t>0,
\\
0 & {\rm if} \quad t\leq 0. \end{array}\right.
\een
There are some basic estimates for $\Gamma_1$ (see, for example, Proposition 2 in [2])
\ben\label{eq:Gamma 1 bound}
\Gamma_1(z,\zeta)\leq C ||\zeta^{-1}\circ z||^{-Q},
\een
and
\ben\label{eq:gradient Gamma 1}
|\ptl_{\xi_i}\,\Gamma_1(z,\zeta)|\leq C ||\zeta^{-1}\circ
z||^{-Q-1},
\een
where $i=1,\cdots,m_0$, for all $z,\zeta\in \mathbb{R}^N\times(0,T]$.

To estimate the operator $\Gamma_1$, we choose to approximate it by considering the homogeneous operator (The $*$ terms of $B$ in $({\bf H_2})$ vanish).
Let $Y_0=\langle x,\,B_0 D\rangle-\partial_ t$, where $B_0 \,$ has the
form
$$\left(
\begin{array}{cccccc}
0 & {B_1} & 0 & \cdots & 0 \\
0 & 0  & {B_2} & \cdots & 0 \\
\vdots & \vdots & \vdots & \ddots &\vdots \\
0 & 0 & 0 & \cdots & {B_d} \\
0 & 0 & 0 & \cdots & 0
\end{array}
\right)
$$
We denote $L_0=div (A_0 D)+Y_0,$ and can define in the same way
$E_0(t)$, ${\mathcal C}_0(t),$ and $\Gamma_0(z,\zeta)$ with respect
to $B_0.$ It is a natural question that what the difference is between ${\mathcal C}_0(t)$ and ${\mathcal C}(t),$ since $B_0$ and $B$ in $({\bf H_2})$ are different. Here we recall that ${\mathcal C}_0(t) (t>0)$ (see[7])
satisfies
\ben\label{eq:bound of C0}
{\mathcal C}_0(t)=\tilde{D}_{t^{\fr 1 2}}{\mathcal C}_0(1)\tilde{D}_{t^{\fr 1 2}}.
\een
Moreover,  the following lemma is obtained by Lanconelli and Polidoro (see Lemma 3.3 in
\cite{LanP}), which states that ${\mathcal C}_0(t)$ and ${\mathcal C}(t)$ are equivalent when $t$ is sufficiently small.

\begin{lemma}\label{lem:C C0} In addition to the above assumptions,  for every given
$T>0$, there exist positive constants $C_T$ and $C'_T$ such that
\ben\label{eq:C t bound}
\langle {\mathcal C}_0(t)x,x\rangle(1-C_T\,t)\leq \langle {\mathcal
C}(t)x,x\rangle \leq \langle {\mathcal
C}_0(t)x,x\rangle(1+C_T\,t),\een
\ben\label{eq:C t -1 bound}
\langle {\mathcal C}_0^{-1}(t)x,x\rangle(1-C_T\,t)\leq \langle
{\mathcal C}^{-1}(t)x,x\rangle \leq \langle {\mathcal
C}_0^{-1}(t)x,x\rangle (1+C_T\,t),\een
\ben\label{eq:det C t bound}
C_T^{'-1}t^Q(1-C_T\,t)\leq {\rm det} {\mathcal C}(t)\leq
C'_T\,t^Q(1+C_T\,t),\een
for every $(x,t)\in \mathbb{R}^N\times (0,T]\,$ and
$t<\,\fr{1}{C_T}.$ \\

\end{lemma}

We need the following classical potential estimates (cf. (1.11) in \cite{Fo}) in our proof.

\begin{lemma}
Let $(\mathbb{R}^{N+1},\circ)$ is a homogeneous Lie group of homogeneous
dimension $Q+2$, $\a \in (0, \frac{Q+2}{p})$ with $p>1$ and $G \in C(\mathbb{R}^{N+1}\setminus
\{0\})$ be a $\delta_{t}$-homogeneous function of degree $\alpha-Q-2$, that is
\beno
G(\delta_t(z))=\delta^{\alpha-Q-2}G(z),\quad z\neq 0.
\eeno
If $f \in L^p(\mathbb{R}^{N+1})$, then
$$
G_f(z)\equiv \int_{\mathbb{R}^{N+1}} G(\zeta ^{-1}\circ z)f(\zeta)d\zeta,
$$
is defined almost everywhere and there exists a constant
$C=C(Q,p)$ such that
$$
||G_f||_{L^q(\mathbb{R}^{N+1})}\leq C \max_{||z||=1} |G(z)|\quad
||f||_{L^p(\mathbb{R}^{N+1})},
$$
where $q$ is defined by
$$
\fr 1q =\fr 1p-\fr{\a}{Q+2}.
$$
\end{lemma}

We also need the following  potential estimates in \cite{CPP}.
\begin{cor}\label{cor:potential estimate}
 Let $f\in L^p(\mathbb{R}^{N+1})$ with $1<p<\infty$.  Recall the definitions in  \cite{CPP}
$$
\Gamma_1(f)(z)=\int_{\mathbb{R}^{N+1}}\Gamma_1(z,\zeta)f(\zeta) d\zeta,
\qquad \forall z\in \mathbb{R}^{N+1},
$$
and
$$
\Gamma_1(D_{m_0}f)(z)=-\int_{\mathbb{R}^{N+1}}D_{m_0}^{(\zeta)}\Gamma_1(z,\zeta)f(\zeta)
d\zeta, \qquad \forall z\in \mathbb{R}^{N+1},
$$
then exists a positive constant $C=C(Q,T,B)$ such that
$$
\|\Gamma_1(f)\|_{L^{q}(S_T)}\leq C\|f\|_{L^p(S_T)},\quad
\fr 1q =\fr 1p-\fr{2}{Q+2}
$$
and
$$
\|\Gamma_1(D_{m_0}f)\|_{L^{q}(S_T)}\leq C\|f\|_{L^p(S_T)},
\quad
\fr 1q =\fr 1p-\fr{1}{Q+2}.
$$
\end{cor}

\setcounter{equation}{0}
\section{Some technical lemmas}

To obtain local estimates of solutions of the equation (\ref{eq:general equ}), for
instance, at point $(x_0,t_0)$, we may consider the estimates at a
ball centered at $0$, since the equation (\ref{eq:general equ}) is invariant
under the left group translation when $a_{ij}$ is constant. By
exploring new  weak ${\rm Poincar\acute{e}}$ type inequality with lower order terms, we prove the
following lemma, Lemma 3.8, which is essential in the oscillation estimates
in Kruzhkov's approaches in parabolic case.

For convenience, let
$x'=(x_1,\cdots,x_{m_0})$ and $x=(x', \overline x)$. Consider the
estimates in the following cube, instead of ${\cal B}^-_r$,
$$
{\cal C}_r^{-}=\{(x,t): \quad-r^2\leq t \leq 0, |x'|\leq r,
|x_{m_0+1}|\leq (\la  r)^{3}, \cdots, |x_N|\leq (\la
r)^{2d+1}\},
$$
where $\lambda$ is a constant satisfying $\lambda\geq \max\{16 N^2,|B|\}$.
Let
$$
K_r=\{x':\quad |x'|\leq r \},
$$
$$
S_r=\{ \overline x:\quad|x_{m_0+1}|\leq (\la  r)^{3},
\cdots, |x_N|\leq (\la  r)^{2d+1}\}.
$$
Moreover, assuming that $0<\a, \beta<1$ are constants, for fixed $t$ and $h$ we assume
$$
{\cal N}_{t,h}\doteq\{(x',\overline x)\in K_{\beta r}\times S_{\beta
r}:\quad u(\cdot,t) \geq h\}.
$$

In the following discussions, we still use the notation
${\cal B}^-_r$ instead of ${\cal C}_r^-$, since there are equivalent, and
we always assume $r \ll 1$. Moreover, all constants
depending on $m_0$, $d$, $N$ or $Q$ will be denoted as $C(B)$.

\begin{lemma}\label{lem:weak subsolution w}
Suppose that $u(x,t)\geq 0$ is a weak solution of equation (\ref{eq:general equ}) in
$\Omega$. Let $w=G(\gamma u+h)$(see the definition of $G$ in Lemma \ref{lem:G function}) with $\gamma>0$ and $0<h\leq \frac14$.  There holds the following inequality
\ben\label{eq:weak subsolution w}
\int_{\Omega} (Dw)^T AD\eta+\eta (Dw)^T ADw-\eta Yw dxdt\leq \int_{\Omega} \left[-\eta b'\cdot D_{m_0}w+\eta \frac{\gamma|cu+f|}{h}\right]dxdt
\een
for any $0\leq \eta\in C_0^{\infty}(\Omega)$.
\end{lemma}

{\it Proof.} The argument follows from the properties of G-function as in \cite{Moser61,Kru64,Gu}, and more details we refer to  Lemma \ref{lem:G function} in the Appendix.

 Let $\phi=\gamma G'(\gamma u+h)\eta$. By (\ref{eq:general equ}) and a smoothing argument, we obatain
\beno
\int_{\Omega} \phi Yu-(Du)^T AD\phi =\int_{\Omega} \phi(b'\cdot D_{m_0}u+ cu+f).
\eeno
Let $w=G(\gamma u+h)$, then $D_{x_i}w=\gamma G'(\gamma u+h)D_{x_i} u$ for $i=1,\cdots,N$ and $D_{t}w=\gamma G'(\gamma u+h)D_{t} u$.
Thus, we get
\beno
\int_{\Omega} \eta Yw-(Dw)^T AD\eta-\eta^2\gamma G''(\gamma u+h)(Du)^T ADu=\gamma\int_{\Omega} G'(\gamma u+h)\eta(b'\cdot D_{m_0}u+ cu+f).
\eeno
By using Lemma \ref{lem:G function}, we get
\beno
G''(\gamma u+h)\geq G'(\gamma u+h)^2,\quad |G'(\gamma u+h)|\leq \frac{1}{h},
\eeno
and we derive that
\beno
\int_{\Omega} (Dw)^T AD\eta+\eta (Dw)^T ADw-\eta Yw\leq \int_{\Omega} -\eta b'\cdot D_{m_0}w+\eta \frac{\gamma|cu+f|}{h}.
\eeno
The proof is complete.
$\diamondsuit$

The following lemma contains an estimate similar to that in \cite{Wang-Zhang2009}.

\begin{lemma}\label{lem:level set}
There exist constants  $\a=\alpha(B)$, $\beta=\beta(B)$, $r_1=r_1(\lambda,B)\leq 1$ and
\beno
h_1=h_1(B,\la,\|b'\|_{L^{Q+2}({\cal B}^-_{1})},\|c\|_{L^q({\cal B}^-_{1})}, \|f\|_{L^q({\cal B}^-_1)}),
\eeno
where $q$ is defined in (\ref{eq:b' condi}), such that for any $h\leq h_1$ and $r^{2-\frac{Q+2}{q}}\leq \min\{r_1^{2-\frac{Q+2}{q}}, h^{\frac98}\}$ we have the following conclusion. Suppose that $u(x,t)\geq 0$ is a solution of equation (\ref{eq:general equ}) in
${\cal B}^-_r$ centered at $0$ and
$$
{\rm meas}\{(x,t)\in {\cal B}^-_r: \quad u \geq 1\} \geq \fr 1 2{\rm  meas }({\cal
B}^-_r).
$$
Then for almost all
$t\in (-\a r^2,0)$, we have
$$
{\rm  meas }\{{\cal N}_{t,h}\} \geq \fr {1}{11}{\rm  meas }\{ K_{\beta r}\times
S_{\beta r}\}.
$$
\end{lemma}
{\it Proof:} Let $\eta(x')$ be a smooth cut-off function so that
\begin{align*}\,\, \left\{
\begin{aligned}
&\eta(x')=1,\quad \hbox {for} \quad |x'|< \beta r,\\
&\eta(x')=0,\quad \hbox {for} \quad |x'|\geq r.
\end{aligned}
\right. \end{align*}
Moreover, $0\leq\eta \leq 1$ and $|D_{m_0} \eta|\leq \fr
{2m_0}{(1-\beta)r}$. Moreover, we let
$$
v=G({u+h^{\fr 9 8}}),
$$
where $h$ is a constant, $0<h<\frac14$, to be decided. Then by Lemma \ref{lem:weak subsolution w} $v$
satisfies
\beno
\int_{{\cal B}^-_r}(Dv)^T AD\psi-\psi Yv+\psi(Dv)^T ADv\leq \int_{{\cal B}^-_r}\psi(-b'\cdot D_{m_0}v+\frac{|cu+f|}{h^{\fr 9 8}}),
\eeno
where $0\leq \psi\in C_0^{\infty}({\cal B}^-_r)$.

Replacing $\psi$ by $\eta^2(x')$ into the above inequality and integrating by parts on
$K_r\times S_{\beta r}\times(\tau,t)$, by a smooth argument to the $t$ or $\bar{x}$ direction we have
\ben\label{eq:energy estimate}
&& \int_{K_{\beta r}}\int_{S_{\beta r}} v(x',\overline x,t)d
\overline x dx' +\fr {1}{2\la}\int_\tau^t \int_{K_{
r}}\int_{S_{\beta r}}\eta^2 \,
|D_{m_0}v|^2d \overline x dx'dt\nonumber\\
&
\leq& \fr {C} {\beta^{Q}(1-\beta)^2}mes(S_{\beta r})mes(K_{\beta
r})+\int_\tau^t
\int_{K_{r}}\int_{S_{\beta r}} \eta^2 x^TBDv d \overline xdx' dt \nonumber\\
&&+\int_{K_{r}}\int_{S_{\beta r}} v(x',\overline x,\tau)d
\overline x dx'+Cr^Q\|b'\|^2_{L^{Q+2}}+Ch^{-\frac98}r^Qr^{2-\frac{Q+2}{q}}\|cu+f\|_q,\quad a.e. \quad\tau, t\in(-r^2,0),\nonumber\\
\een
where $C$ only depends on $\la$ and $B$. Let
$$
I_B\equiv\int_{K_{r}}\int_{S_{\beta r}} \eta^2 \sum_{i,j=1}^N
x_ib_{ij}\ptl_{x_j}v d \overline x dx'= I_{B_1}+I_{B_2},
$$
where
$$
I_{B_1}=\int_{K_{r}}\int_{S_{\beta r}} \eta^2 \sum_{i=1}^N
\sum_{j=1}^{m_0} x_ib_{ij}\ptl_{x_j}v d \overline x dx',
$$
$$
I_{B_2}=\int_{K_{r}}\int_{S_{\beta r}} \eta^2 \sum_{i=1}^N
\sum_{j=m_0+1}^{N} x_ib_{ij}\ptl_{x_j}v d \overline x dx'.
$$
By straightforward computations, we see that
\beno
|I_{B_1}|&\leq&\int_{K_{r}}\int_{S_{\beta r}}
\varepsilon\eta^2|D_{m_0}v|^2+C_\varepsilon\eta^2\sum_{j=1}^{m_0}\sum_{i=1}^N
|x_i b_{ij}|^2 d \overline x dx' \\
&\leq&\int_{K_{r}}\int_{S_{\beta r}}
\varepsilon\eta^2|D_{m_0}v|^2d\overline x
dx'+C(\varepsilon,B,\lambda)\beta^{-Q}|K_{\beta r}||S_{\beta
r}|,\eeno
and noting that $\lambda>|B|$
$$
\ba{llllllll} |I_{B_2}| &\leq & |\int_{K_{r}}\int_{S_{\beta r}}
\eta^2 \sum_{i=1}^N \sum_{j=m_0+1}^{N} x_i b_{ij}\ptl_{x_j}v
d\overline x dx'| \\ \\& \leq & |\int_{K_{r}}\int_{S_{\beta r}}
{-}\eta^2\sum_{i=1}^N\sum_{j>m_0}\delta_{ij}b_{ij}vd \overline x dx'|\\
\\&&+|\sum_{i=1}^N\sum_{j>m_0}\int_{K_{r}}\int_{\partial_j S_{\beta r}}
\eta^2x_{i}b_{ij}v d\overline {x_j} dx'| \\
\\& \leq & \lambda N\beta^{-Q}|K_{\beta r}||S_{\beta r}| \ln
(h^{-\fr 9 8})\\ \\&& +\la\sum_{i=1}^N\sum_{j>m_0}\fr{(\lambda
r)^{\alpha_i}}{(\lambda r)^{\alpha_j}} \beta^{-2Q}|K_{\beta
r}||S_{\beta r}| \ln (h^{-\fr 9 8}), \ea
$$
where $\overline {x_j}=(x_{m_0+1}, \dots, x_{j-1}, x_{j+1}, \dots,
x_N)$ and $\partial_j S_{\beta r}$ denotes the boundary of $S_{\beta r}$ with respect to the $j$-th component. When $\a_i\geq\a_j$, we have
$$
\int_\tau^t |I_{B_2}| \leq (\lambda N \,r^2+\la r^2
N^2)\beta^{-2Q}|K_{\beta r}||S_{\beta r}| \ln (h^{-\fr 9 8}),
$$
or $i<j$, thus $\a_j=\a_i+2$ by the property of $B$, then
$$
\int_\tau^t |I_{B_2}| \leq (\lambda N \,r^2+\la^{-1}
N^{2})\beta^{-2Q}|K_{\beta r}||S_{\beta r}| \ln (h^{-\fr 9 8}).
$$
By $\la>16 N^2$, we choose $r_1$ small enough, such that for any $r\leq r_1$
$$\lambda N \,r^2+2\la r^2N^2+\lambda^{-1}N^{2}<\fr{1}{8},$$ thus
$$
\int_\tau^t |I_{B_2}|\leq \fr{1}{4}\beta^{-2Q}|K_{\beta r}||S_{\beta
r}| \ln (h^{-\fr 9 8}).
$$
By integration to time for $I_B$, we have
\ben\label{eq:estimate of B}
&&\hspace*{-10pt}\int_\tau^t\int_{K_{r}}\int_{S_{\beta r}}
\eta^2 x^T B Dv d \overline x dx'dt \nonumber\\
& \leq& \fr {1}{4}
{\beta}^{-2Q}
\ln(h^{-\fr 9 8}) mes(S_{\beta r})mes(K_{\beta r})\nonumber\\
&&+\int_\tau^t\int_{K_{r}}\int_{S_{\beta r}}
\varepsilon\eta^2|D_{m_0}v|^2+C(\varepsilon,B,\lambda)\beta^{-Q}|K_{\beta
r}||S_{\beta r}|.\een

We shall estimate the measure of the set ${\cal N}_{t,h}$. Let
$$
\mu(t)=mes\{(x',\overline x)|\quad x'\in K_r,\, \overline x \in
S_{r}, \, u(\cdot,t)\geq 1\}.
$$
By our assumption, for $0<\a< \fr 12$
$$
\fr 12 r^2 mes(S_{r})mes(K_{r})\leq \int_{-r^2}^0
\mu(t)dt=\int_{-r^2}^{-\a r^2}\mu(t)dt+\int_{-\a r^2}^{0}\mu(t)dt,
$$
that is
$$
\int_{-r^2}^{-\a r^2}\mu(t)dt\geq (\fr 12-\a)r^2
mes(S_{r})mes(K_{r}),
$$
then there exists a $\tau \in (-r^2,-\a r^2)$, such that
$$
\mu(\tau)\geq (\fr 12-\a)(1-\a)^{-1}
mes(S_{r})mes(K_{r}),
$$
we have by noticing $v=0$ when $u\geq 1,$
$$
\int_{K_{r}}\int_{S_{\beta r}} v(x',\overline x,\tau)d \overline x
dx'\leq \fr 12(1-\a)^{-1}mes(S_{r})mes(K_{r})\ln(h^{-\fr 9
8}).
$$
Now we choose $\varepsilon={\fr{1}{2\la}}\,$ and $\a$  (near zero)
and $\beta$ (near one), so that
\ben\label{eq:choose of alpha beta}
\fr{1}{4\beta^{2Q}}+\fr{1}{2\beta ^{2Q}(1-\a)}\leq \fr 4
5.
\een
Note that the last two terms of (\ref{eq:energy estimate}) and the last term in (\ref{eq:estimate of B})
can be controlled by
$$C(B,\la,\|b'\|_{L^{Q+2}},\|c\|_{L^q}, \|f\|_{L^q})(1-\beta)^{-2}\beta^{-Q}|K_{\beta
r}| | S_{\beta r}|$$
by choosing $r^{2-\frac{Q+2}{q}}\leq h^{\fr98}$.

Combining (\ref{eq:energy estimate})-(\ref{eq:choose of alpha beta}), we deduce
\beno
\int_{K_{\beta r}}\int_{S_{\beta r}} v(x',\overline x,t)d \overline
x dx'\leq [2C(1-\beta)^{-2}\beta^{-Q} +\fr 45\ln(h^{-\fr 9
8})]mes(K_{\beta r}\times S_{\beta r}).\eeno
When $(x', \bar{x})\notin {\cal N}_{t,h},$ $u\geq h$, we have
$$\ln(\fr 1 {2h})\leq \ln^+(\fr{1}{h+h^{\fr 9 8}})\leq v,$$
then
$$\ln(\fr 1
{2h})mes(K_{\beta r}\times S_{\beta r}\setminus {\cal N}_{t,h})\leq
\int_{K_{\beta r}}\int_{S_{\beta r}} v(x',\overline x,t)d \overline
x dx'.$$
Since
$$
\fr{C+{\fr 45}\ln(h^{-\fr 98})}{\ln(h^{-1})}\lra \fr
9{10},\qquad\hbox{as} \quad h\ra 0,
$$
then there exists constant $h_1$ such that for $0<h<h_1$ and $t
\in(-\a r^2,0)$
$$
mes(K_{\beta r}\times S_{\beta r}\setminus {\cal N}_{t,h})\leq \fr
{10}{11}mes(K_{\beta r}\times S_{\beta r}).
$$
Hence the proof is complete.$\diamondsuit$

Let $\chi(s)$ be a smooth function given by
$$\ba{ll}
\chi(s)=1 \qquad if \quad s\leq {\theta^{\fr 1 {Q}}} r,\\
\chi(s)=0 \qquad if \quad s> r, \ea
$$
where ${\theta^{\fr 1 {Q}}}<\fr {1}{2}$ is a constant. Moreover, we
assume that
$$
0\leq -\chi'(s) \leq \fr{2}{(1 -{\theta^{\fr 1 {Q}}})r},
$$
and $\chi'(s)<0$, if ${\theta^{\fr 1 {Q}}} r<s< r$. Also for any
$\beta_1, \beta_2,$ with $\theta^{\fr 1 {Q}}<\beta_1<\beta_2<1,$ we
have $$|\chi'(s)|\geq C(\beta_1,\beta_2)>0,$$ if $\beta_1r\leq s\leq
\beta_2r.$

For $x\in \mathbb{R}^N,$ $t<0$, we set $${\mathcal Q} =\{(x',\bar{x},t)|
-r^2\leq t< 0,\, x'\in K_{\fr r \theta},
\,|x_j|\leq\fr{r^{\a_j}}{\theta}, j=m_0+1,\cdots, N\},$$
$$
\phi_0(x,t)=\chi\left(\left[\theta^2\sum_{i=m_0+1}^{N} \fr
{x_i^2}{r^{2\a_i-Q}}-C_1 t r^{Q-2} \right]^{\fr {1} {Q}}\right),
$$
$$
\phi_1(x,t)=\chi(\theta |x'|),
$$
\ben\label{eq:def of phi}
\phi(x,t)=\phi_0(x,t)\phi_1(x,t),
\een
where $C_1>1$ is chosen so that
$$
\ba{lllll}  C_1 r^{Q-2}&\geq \theta^2|\sum_{i=1}^{N}\sum_{j>m_0}2x_i
b_{ij}x_j r^{Q-2\alpha_j}|, \ea
$$
for all $z\in {\mathcal Q}$.

In the following discussion, $a\approx b$ means
$$C(B,\la)^{-1}a\leq b\leq C(B,\la)a.$$

%Applying (\ref{eq:C t -1 bound}) and (\ref{eq:bound of C0}),
%\begin{eqnarray*}
%\langle {\mathcal C}^{-1}(|t|)e^{t B^T}x,e^{t B^T}x\rangle &\approx
%&\langle {\mathcal C}_0^{-1}(|t|)e^{t B^T}x,e^{t B^T}x\rangle\\ &=&
%\langle {\mathcal C}_0^{-1}(1)D_{|t|^{-\fr 1
%2}}e^{t\,B^T}x,D_{|t|^{-\fr 1 2}}e^{t\,B^T}x\rangle\\ &\approx
%&\|e^{\tilde{B}}D_{|t|^{-\fr 1 2}}x\|\approx |D_{|t|^{-\fr 1 2}}x|^2
%\end{eqnarray*}
%where $D_{|t|^{-\fr 1 2}}B^T=|t|^{-1}\tilde{B}D_{|t|^{-\fr 1 2}},\,$
%$D_{|t|^{-\fr 1 2}}e^{tB^T}=e^{\tilde{B}} \, D_{|t|^{-\fr 1 2}}\,$
%and $\tilde{B}$\, has the form
%$$\left(
%\begin{array}{ccccc}
%|t|B_{0,0}^T  &|t|^2 B_{1,0}^T & \cdots & \cdots & |t|^{d+1}B_{d,0}^T  \\
%B_1^T &  |t|B_{1,1}^T   &\cdots &\cdots & |t|^{d}B_{d,1}^T \\
%0 &B_2^T& \ddots &\cdots &\vdots \\
%\vdots & \ddots &\ddots &\ddots & \vdots \\
%0 &  \cdots  & 0 &  B_d^T &|t|B_{d,d}^T \\
%\end{array}
%\right)
%$$
%B is given by
%$$\left(
%\begin{array}{ccccc}
%B_{0,0}&B_1&0&\cdots&0\\
%B_{1,0}&B_{1,1}&B_2&\cdots&0\\
%\vdots &\vdots &\vdots &\ddots &\vdots \\
%B_{d-1,0}&B_{d-1,1}&B_{d-1,2}&\cdots&B_d\\
%B_{d,0}&B_{d,1}&B_{d,2}&\cdots&B_{d,d}\\
%\end{array}
%\right)
%$$
%then we obtain (a).
\begin{rem}[c.f. Remark 3.1 in  \cite{Wang-Zhang2009}] By the definition of $\phi$ and the above arguments,
it is easy to check that, for $\theta$,  $r$ small enough and $t\leq 0$\\
(1) $\phi(z)\equiv 1,$ in ${\cal B}^-_{\theta r}$,\\
(2) $\rm{supp}\phi\bigcap\{(x,t);t\leq 0\}\subset {\mathcal Q}$,\\
(3) there exists $\a_1>0,$ which depends on $C_1,$ such that
$$\{(x,t)| -\a_1r^2\leq t < 0, x'\in K_r, \bar{x}\in S_{\beta
r}\}\subseteq \rm{supp}\phi, $$ (4) $0<\phi_0(z)<1,$ for $z\in
\{(x,t)| -\a_1r^2\leq t \leq -\theta r^2, x'\in K_r, \bar{x}\in
S_{\beta r}\}$. Note that $\alpha_1>\theta$ if $theta$ is chosen to be small enough.
\end{rem}

Using Lemma \ref{lem:C C0} and the properties $\phi$, we have the following lemma (Note that $B$ and $Y$ are as in \cite{Wang-Zhang2009}).
 \begin{lemma}[c.f. (d.1) and Lemma 3.2 in \cite{Wang-Zhang2009}, p1598-1601]\label{lem:Y phi} Under the above notations, we have\\
(a) For $t<0$, $|t|$ is small
enough, then we have
$$
\hspace*{12pt}\langle {\mathcal C}^{-1}(|t|)e^{t B^T}x,e^{t
B^T}x\rangle\approx |\tilde{D}_{|t|^{-\fr 1 2}}x|^2,$$
where C depends on $B$ and $\la$.\\
(b) $Y \phi_0(z)\leq 0, \quad \rm{for}\quad z\in {\mathcal Q}$.
\end{lemma}
%{\it Proof:}Let
%$
%[\theta^2\sum_{i=m_0+1}^{N} \fr {x_i^2}{r^{2\a_i-Q}}-C_1 t r^{Q-2}]
%$
%be denoted by $[\cdots]$. Then
%$$\ba{llllllllll} Y \phi_0 &=
%\chi'([\cdots]^{\fr{1}{Q}})\fr{1}{Q}[\cdots]^{\fr{1}{Q}-1} [C_1
%r^{Q-2} +\theta^2\sum_{i=1}^{N}\sum_{j>m_0}(2x_i b_{ij}x_j
%r^{Q-2\alpha_j})]
% \ea
%$$
%We choose $C_1>1,$ such that
%$$
%\ba{lll}  C_1 r^{Q-2}\geq \theta^2|\sum_{i=1}^{N}\sum_{j>m_0}2x_i
%b_{ij}x_j r^{Q-2\alpha_j}|, \ea
%$$
%
%For the term $x_i b_{ij}x_j r^{Q-2\alpha_j}$, either $\a_i\geq \a_j$
%or $\a_j=\a_i+2$, we obtain
%$$ \theta^2|\sum_{i=1}^{N}\sum_{j>m_0}2x_i b_{ij}x_j
%r^{2Q-2\alpha_j}| \leq C(B,\la) r^{Q-2}.
%$$
%Thus $C_1(B,\la)$ is well defined,
%then $Y \phi_0(z)\leq 0$ ($z\in {\mathcal Q}$) holds. \\

Let $
w=G(\frac{u}{h}+h^{\fr 18}).
$ Then we have the following ${\rm Poincar\acute{e}}$'s type
inequality.
\begin{lemma}[Weak Poincar\'{e} inequality]\label{lem:weak poincare}
Let $u$ be a non-negative weak solution of (\ref{eq:general equ}) in ${\cal
B}_1^-$ and $
w=G(\frac{u}{h}+h^{\fr 18}).
$ Then there exists a constant $C=C(B,\la)$ such that for $r<\theta<1$
\ben\label{eq:weak poincare}
&&\int_{{\cal B}^-_{\theta r}}(w(z)-I_0)_+^2
\leq C\theta^2
r^2\int_{{\cal B}^-_{\fr r {\theta}}}|D_{m_0}w|^2\nonumber\\
&&+C(B,\la) h^{-\fr94}|\fr{r}{\theta}|^{Q+2}|\fr{r}{\theta}|^{\fr{8}{Q+2}-\frac4q}\left(||c
||^2_{L^{q}({\cal B}^-_{{\fr r {\theta}}})}||u
||^2_{L^{\infty}({\cal B}^-_{{\fr r {\theta}}})}+||f
||^2_{L^{q}({\cal B}^-_{{\fr r {\theta}}})}\right),
\een
where $I_0$ is given by
\ben\label{eq:I0}
I_0={\rm max}_{{\cal B}^-_{\theta r}}[I_1(z)+C_2(z)],
\een
and
\ben\label{eq:I1 C2}
I_1(z)&=&\int_{{\cal B}^-_{\fr r {\theta}}} [\langle
{\phi}_1A_0D{\phi}_0,D\Gamma_1(z,\cdot)\rangle
w-\Gamma_1(z,\cdot)wY\phi](\zeta)d\zeta,\nonumber\\
C_2(z)&=&\int_{{\cal B}^-_{\fr r {\theta}}} [\langle
{\phi}_0A_0D{\phi}_1,D\Gamma_1(z,\cdot)\rangle w](\zeta)d\zeta,
\een
where $\Gamma_1$ is the fundamental solution,  and $\phi$ is given
by (\ref{eq:def of phi}).
\end{lemma}
 {\it Proof:} We represent $w$ in terms of the fundamental
solution of $\Gamma_1$. For $z \in {\cal B}^-_{\theta r}$, we have
\ben\label{eq:w equality}
w(z)&=&\int_{{\cal B}^-_{\fr r {\theta}}}  [\langle
A_0D(w\phi),D\Gamma_1(z,\cdot)\rangle
-\Gamma_1(z,\cdot)Y(w\phi)](\zeta)d\zeta \nonumber\\
&=&
I_1(z)+I_2(z)+I_3(z)+C_2(z),\een
where $I_1(z)$ and $C_2(z)$ are given by (\ref{eq:I1 C2}) and
\beno
&&I_2(z)=\int_{{\cal B}^-_{\fr r {\theta}}} [\langle
(A_0-A)Dw,D\Gamma_1(z,\cdot)\rangle\phi-\Gamma_1(z,\cdot)\langle
ADw,D\phi\rangle](\zeta)d\zeta\\
&&+\int_{{\cal B}^-_{\fr r {\theta}}} [-\Gamma_1(z,\cdot)\phi b'\cdot D_{m_0}w+\Gamma_1(z,\cdot)\phi h^{-\fr98}(|cu|+|f|)](\zeta)d\zeta\\
&\equiv&I_{21}+\cdots+I_{25},
\eeno
\beno
I_3(z)=\int_{{\cal B}^-_{\fr r {\theta}}} [\langle
ADw,D(\Gamma_1(z,\cdot)\phi)\rangle-\Gamma_1(z,\cdot)\phi
Yw+\Gamma_1(z,\cdot)\phi b'\cdot D_{m_0}w-\Gamma_1(z,\cdot)\phi\frac{|cu|+|f|}{h^{\fr98}}](\zeta)d\zeta.
\eeno
From our assumption, $w$ satisfies (\ref{eq:weak subsolution w}), and $\phi(\zeta)\Gamma_1(z,\cdot)$
is a test function of this semi-cylinder. In fact, we let
$$
\tilde{\chi}(\tau)=\left\{
\begin{array}{lll} 1\quad &\tau\leq
0,\\1-n\tau\quad &0\leq \tau \leq 1/n,\\ 0\quad &\tau\geq
1/n.\end{array}\right.$$
Then
$\tilde{\chi}(\tau)\phi\Gamma_1(z,\cdot)$ can be a test function
(see [2]). Let $n\rightarrow \infty$, we obtain
$\phi\Gamma_1(z,\cdot)$ as a legitimate test function, and
$I_3(z)\leq 0$. Then in ${\cal B}^-_{\theta r}$,
$$
0\leq (w(z)-I_0)_+\leq |I_2(z)|\leq |I_{21}|+\cdots+|I_{25}|.
$$
By Corollary \ref{cor:potential estimate}, we have
\beno
||I_{21}||_{L^2({\cal B}^-_{\theta r})}\leq C(\la)\theta
r||I_{21}||_{L^{2+\fr 4 Q}({\cal B}^-_{\theta r})}\leq C(B, \la)
\theta r||D_{m_0}w||_{L^2({\cal B}^-_{{\fr r
{\theta}}})}.\eeno
Similarly for $I_{22},$
$$
||I_{22}||_{L^2({\cal B}^-_{\theta r})}\leq |{\cal B}^-_{\theta
r}|^{\fr 12-\fr{Q-2}{2Q+4}} ||I_{22}||_{L^{2\tilde{k}}({\cal
B}^-_{\theta r})}\leq C(B,\la) \theta^2 r^2||D_{m_0}w
D_{m_0}\phi||_{L^2({\cal B}^-_{{\fr r {\theta}}})},
$$
where $|D_{m_0}\phi|=|\phi_0
D_{m_0}\phi_1|=|\phi_0\chi'(\theta|\xi'|)\theta D_{m_0}(|\xi'|)|\leq
\fr{\theta}{r}, $ thus $$||I_{22}||_{L^2({\cal B}^-_{\theta r})}\leq
C(B,\la)
\theta^2 r||D_{m_0}w||_{L^2({\cal B}^-_{{\fr r {\theta}}})}.$$

For $I_{23}$, we have
\beno
||I_{23}||_{L^2({\cal B}^-_{\theta r})}&\leq& |{\cal B}^-_{\theta
r}|^{\fr{Q}{Q+2}} ||I_{23}||_{L^{\fr{2Q+4}{Q}}({\cal
B}^-_{\theta r})}\\
&\leq& C(B,\la) \theta^2 r^2||b'\cdot D_{m_0}w
||_{L^{\fr{Q+4}{2Q+4}}({\cal B}^-_{{\fr r {\theta}}})}\\
&\leq& C(B,\la) \theta^2 r^2|| b'
||_{L^{Q+2}({\cal B}^-_{{\fr r {\theta}}})} || D_{m_0}w
||_{L^{2}({\cal B}^-_{{\fr r {\theta}}})}
\eeno

For  $I_{24}$, we have
\beno
||I_{24}||_{L^2({\cal B}^-_{\theta r})}
&\leq& C(B,\la) h^{-\fr98}||cu
||_{L^{\fr{2Q+4}{Q+6}}({\cal B}^-_{{\fr r {\theta}}})}\\
&\leq& C(B,\la) h^{-\fr98}|\fr{r}{\theta}|^{\fr{Q+2}{2}}|\fr{r}{\theta}|^{\fr{4}{Q+2}-\frac2q}||c
||_{L^{q}({\cal B}^-_{{\fr r {\theta}}})}||u
||_{L^{\infty}({\cal B}^-_{{\fr r {\theta}}})}.
\eeno

For  $I_{25}$, we have
\beno
||I_{25}||_{L^2({\cal B}^-_{\theta r})}
&\leq& C(B,\la) h^{-\fr98}||f
||_{L^{\fr{2Q+4}{Q+6}}({\cal B}^-_{{\fr r {\theta}}})}\\
&\leq& C(B,\la) h^{-\fr98}|\fr{r}{\theta}|^{\fr{Q+2}{2}}|\fr{r}{\theta}|^{\fr{4}{Q+2}-\frac2q}||f
||_{L^{q}({\cal B}^-_{{\fr r {\theta}}})}.
\eeno
Then we prove our lemma.
   $\diamondsuit$

Next, we'll sketch the proof of the weak Sobolev inequality and $L^{\infty}$ bounded estimates as in \cite{CP,CPP}.
In fact, we obtain two types of weak Sobolev equalities, where the representation formula of fundamental solution and potential estimates in Corollary  \ref{cor:potential estimate} are used.

\begin{lemma}[Sobolev estimate]\label{lem:Sobolev estimate}
Under the assumptions $({\bf H_1}-{\bf H_3})$,
 let $u$ be a non-negative weak sub-solution of (\ref{eq:general equ}) in $\Omega$.\\
(i) For $(x_0,t_0)\in \Omega$ and $\overline{{\cal
B}^-_r(x_0,t_0)}\subset \Omega$, there holds
\ben\label{eq:Sobolev estimate}
\|\varphi u\|_{L^{2k}({\cal B}^-_{\rho}(x_0,t_0))}\leq \frac{C}{r-\rho}(\|u\|_{L^2({\cal B}^-_{r}(x_0,t_0))}+ \|D_{m_0}u\|_{L^2({\cal B}^-_{r}(x_0,t_0))})+C\|f\|_{L^{\frac{2Q+4}{Q+4}}({\cal B}^-_{r}(x_0,t_0))}
\een
where $\varphi$ be a cut-off function such that $\varphi=1$ in ${\cal B}^-_{\rho}$, $\frac12 \leq \rho<r\leq 1,$ $k=1+\frac{2}{Q}$, and $C$ depends only on $q,N,\la,Q, \|b'\|_{L^{Q+2}({\cal B}^-_{r}(x_0,t_0))}$ and $\|c\|_{L^{q}({\cal B}^-_{r}(x_0,t_0))}.$\\
(ii) Moreover, let $w=u^p$ with the positive integer $p> 1$, and we have the following similar estimate
\ben\label{eq:Sobolev estimate-w}
\|\varphi w\|_{L^{2k}({\cal B}^-_{\rho}(x_0,t_0))}&\leq& \frac{C}{r-\rho}(p^{\frac{2}{1-\tilde{\beta}}}\|w\|_{L^2({\cal B}^-_{r}(x_0,t_0))}+ \|D_{m_0}w\|_{L^2({\cal B}^-_{r}(x_0,t_0))})\nonumber\\
&&+ p^{-p}r^{p(2-\frac{Q+2}{q})+\frac{Q}{2}}\|f\|_{L^q({\cal B}^-_{r}(x_0,t_0))}^p
\een
where $\tilde{\beta}=\frac{Q+2}{q}-1\in (0,1)$ with $q>(Q+2)/2$.
%(iii) Furthermore, if we assume that $f\in L^{Q+2}_{loc}$, then we have
%\ben\label{eq:Sobolev estimate-ww}
%\|w\varphi\|_{L^{2k}}&\leq& \frac{C}{r-\rho}\left(p^{\frac{1}{1-\beta}}\|w\|_{L^2({\cal B}^-_{r})}+ \|D_{m_0}w\|_{L^2({\cal B}^-_{r})}\right)+r^{(Q+2)/2}\|f\|_{L^{q}({\cal B}^-_{r})}^{p}.
%\een
\end{lemma}

{\bf Proof of Lemma \ref{lem:Sobolev estimate}:}
{\bf Step I: Test function.} Following the method introduced in Lemma 3 in \cite{CPP}, we use $\Gamma_1(z,\cdot)\varphi$ as a test function in the definition of weak solution (\ref{eq:weak subsolution}), which is made by the cut-off at the singularity and dominated convergence theorem.
For example, for the term
$cu\in L^{\frac{2q}{q+2}}$ with $\frac{2q}{q+2}> 1$, since $q>\frac{Q+2}{2}\geq 2$, we have $\int_{\Omega}\Gamma_1(z,\zeta)\varphi(\zeta) c(\zeta)u(\zeta)d\zeta\in L^m$ with $\frac1m=\frac12+\frac1q-\frac{2}{Q+2}$ due to Corollary \ref{cor:potential estimate}, and obviously $m>2$. Hence, we get
\beno
\int_{\Omega}\Gamma_1(z,\cdot)\varphi cu\chi(\frac{\|\zeta^{-1}\circ z\|}{\varepsilon})\rightarrow \int_{\Omega}\Gamma_1(z,\cdot)\varphi cu,
\eeno
for almost every $z\in \mathbb{R}^{N+1}$, where $\chi$ is a smooth function satisfying $\chi(s)=0$ for $s\in [0,1]$ and $\chi(s)=1$ for $s\geq 2$.
For the others, we omitted it. Consequently, we get
\ben\label{eq:weak subsolution-gamma}
\int_{\Omega} \Gamma_1(z,\cdot)\varphi Yu-(Du)^T AD(\varphi\Gamma_1(z,\cdot)) \geq\int_{\Omega} \Gamma_1(z,\cdot)\varphi(b'\cdot D_{m_0}u+ cu+f).
\een

{\bf Step II: Proof of (\ref{eq:Sobolev estimate}).} Let $\varphi$ be a cut-off function such that $\varphi=1$ in ${\cal B}^-_{\rho}$ and $\varphi=0$ outside of ${\cal B}^-_r$; furthermore, $|\partial_t \varphi|+|D\varphi|\leq \frac{C}{r-\rho}.$
We represent $u$ in terms of the fundamental
solution of $\Gamma_1$. For $z \in {\cal B}^-_{\rho}$, we have
\beno
u(z)\varphi(z)&=&\int_{{\cal B}^-_r}  [\langle
A_0D(u\varphi),D\Gamma_1(z,\cdot)\rangle
-\Gamma_1(z,\cdot)Y(u\varphi)](\zeta)d\zeta
\\
&=&
I_1(z)+I_2(z)+I_3(z)+I_4(z),
\eeno
where $I_1(z)-I_4(z)$ are as follows:
\beno
I_1(z)=\int_{{\cal B}^-_r} [\langle
A_0D{\varphi},D\Gamma_1(z,\cdot)\rangle
u-\Gamma_1(z,\cdot)uY\varphi](\zeta)d\zeta,
\eeno
\beno
I_2(z)=\int_{{\cal B}^-_r} [\langle
(A_0-A)Du,D\Gamma_1(z,\cdot)\rangle\varphi-\Gamma_1(z,\cdot)\langle
ADu,D\varphi\rangle](\zeta)d\zeta,
\eeno
\beno
I_3(z)=\int_{{\cal B}^-_r} [\langle
ADu,D(\Gamma_1(z,\cdot)\varphi)\rangle-\Gamma_1(z,\cdot)\varphi
Yu+\Gamma_1(z,\cdot)\varphi (b'\cdot D_{m_0}u+cu+f)](\zeta)d\zeta'
\eeno
and
\beno
I_4(z)=-\int_{{\cal B}^-_r} [\Gamma_1(z,\cdot)\varphi (b'\cdot D_{m_0}u+cu+f)](\zeta)(\zeta)d\zeta=I_{41}(z)+I_{42}(z)+I_{43}(z),
\eeno

Obviously, by (\ref{eq:weak subsolution-gamma}) we have
$
I_3(z)
\leq 0.
$
For the term $I_1$, by  Corollary \ref{cor:potential estimate} and H\"{o}lder inequlity we have
\beno
\|I_1\|_{L^{2k}({\cal B}^-_r)}\leq C\|u D_{m_0}\varphi\|_{L^2(R^{N+1})}+C|{\cal B}^-_r|^{\frac{2}{Q}}\|u Y\varphi\|_{L^2(R^{N+1})}\leq \frac{C}{r-\rho}\|u\|_{L^2({\cal B}^-_{r})}
\eeno
Similarly, for $I_2$ we have
\beno
\|I_2\|_{L^{2k}({\cal B}^-_r)}\leq \frac{C}{r-\rho}\|D_{m_0}u\|_{L^2({\cal B}^-_{r})}
\eeno

Next, we estimate $I_4$,
\beno
\|I_{41}\|_{L^{2k}({\cal B}^-_r)}\leq C\|\varphi b'\cdot D_{m_0}u\|_{L^{\frac{2Q+4}{Q+4}}(R^{N+1})}\leq C\|b'\|_{L^{Q+2}({\cal B}^-_{r})}\|D_{m_0}u\|_{L^2({\cal B}^-_{r})},
\eeno
\beno
\|I_{43}\|_{L^{2k}({\cal B}^-_r)}\leq C\|\varphi f\|_{L^{\frac{2Q+4}{Q+4}}(R^{N+1})}\leq C\|f\|_{L^{\frac{2Q+4}{Q+4}}({\cal B}^-_{r})},
\eeno
and
\beno
\|I_{42}\|_{L^{2k}({\cal B}^-_r)}&\leq& C\|c|\varphi u|^{\tilde{\beta}}|\varphi u|^{1-\tilde{\beta}}\|_{L^{\frac{2Q+4}{Q+4}}(R^{N+1})}\leq C\|c\|_{L^{q}({\cal B}^-_{r})}\|\varphi u\|_{L^{2k}({\cal B}^-_{r})}^{\tilde{\beta}}\|u\varphi\|_{L^2({\cal B}^-_{r})}^{1-\tilde{\beta}}\\
&\leq& \frac12\|\varphi u\|_{L^{2k}({\cal B}^-_{r})}+C\|c\|_{L^{q}({\cal B}^-_{r})}^{\frac{1}{1-\tilde{\beta}}}\|u\varphi\|_{L^2({\cal B}^-_{r})},
\eeno
where $\tilde{\beta}=\frac{Q+2}{q}-1\in (0,1)$.

Concluding the above estimates, we have
\beno
\|u\varphi\|_{L^{2k}}&\leq& \frac{C}{r-\rho}\|u\|_{L^2({\cal B}^-_{r})}+ C(\frac{1}{r-\rho}+ \|b'\|_{L^{Q+2}({\cal B}^-_{r})})\|D_{m_0}u\|_{L^2({\cal B}^-_{r})}\\
&&+C\|f\|_{L^{\frac{2Q+4}{Q+4}}({\cal B}^-_{r})}+C\|c\|_{L^{q}({\cal B}^-_{r})}^{\frac{q}{2q-Q-2}}\|u\varphi\|_{L^2({\cal B}^-_{r})}.
\eeno
Then we complete the proof of  (\ref{eq:Sobolev estimate}).

{\bf Step III: Proof of (\ref{eq:Sobolev estimate-w}).} To prove the inequality (\ref{eq:Sobolev estimate-w}), if
\beno
\|w\|_{L^2({\cal B}^-_{r}(x_0,t_0))}+ \|D_{m_0}w\|_{L^2({\cal B}^-_{r}(x_0,t_0))}<\infty,
 \eeno
at first the function $pu^{p-1}\Gamma_1(z,\cdot)\varphi$ can be seen as a test function of (\ref{eq:general equ}) and there holds
\ben\label{eq:weak subsolution-gamma-}
\int_{\Omega} \Gamma_1(z,\cdot)\varphi Yw-(Dw)^T AD(\varphi\Gamma_1(z,\cdot)) \geq\int_{\Omega} \Gamma_1(z,\cdot)\varphi(b'\cdot D_{m_0}w+ pcw+pu^{p-1}f).
\een
Next, we deal with the terms $pcu^p$ and $pu^{p-1}|f|$ only, and other terms are similar as Step II. Write the last two terms of the righthand of (\ref{eq:weak subsolution-gamma-}) as $I_{42}'$ and $I_{43}'$.

{\bf Case I:  $f\in L^q$ and $p>\frac{q(Q+2)}{2(Q+2-q)}$.} At this moment, we have
\beno
\|I_{42}'\|_{L^{2k}({\cal B}^-_r)}
\leq \frac12\|\varphi w\|_{L^{2k}({\cal B}^-_{r})}+Cp^{\frac{1}{1-\tilde{\beta}}}\|c\|_{L^{q}({\cal B}^-_{r})}^{\frac{1}{1-\tilde{\beta}}}\|w\varphi\|_{L^2({\cal B}^-_{r})},
\eeno
where $\tilde{\beta}=\frac{Q+2}{q}-1\in (0,1)$. And by Corollary 2.1
\beno
\|I'_{43}\|_{L^{2k}({\cal B}^-_r)}&\leq& C\|\varphi pu^{p-1} f\|_{L^{\frac{2Q+4}{Q+4}}(R^{N+1})}\leq Cp\|f\|_{L^{q}({\cal B}^-_{r})}\|\varphi w\|_{L^{2k}({\cal B}^-_{r})}^{\beta'}\|w\varphi\|_{L^2({\cal B}^-_{r})}^{1-\frac1p-\beta'}\\
&\leq& \frac12\|\varphi w\|_{L^{2k}({\cal B}^-_{r})}+\frac{C}{r-\rho}p^{\frac{2}{1-\frac1p-\beta'}}\|w\varphi\|_{L^2({\cal B}^-_{r})}+p^{-p}(r-\rho)^{p(1-\frac1p-\beta')}\|f\|_{L^{q}({\cal B}^-_{r})}^{p}
\eeno
where $\beta'$ satisfies
\beno
\frac1q+\frac{\beta'}{2k}+\frac{1-\frac1p-\beta'}{2}=\frac{Q+4}{2Q+4},\quad k=1+\frac{2}{Q}
\eeno
and hence $\beta'=(\frac1q-\frac{1}{2p})(Q+2)-1\in (0,1)$ with $q>(Q+2)/2$.

It is easy to check that
\beno
\frac{1}{1-\tilde{\beta}}>\frac{1}{1-\frac1p-\beta'},
\eeno
which yields that
\beno
\|w\varphi\|_{L^{2k}({\cal B}^-_{r})}&\leq& \frac{C}{r-\rho}p^{\frac{2}{1-\tilde{\beta}}}\|w\|_{L^2({\cal B}^-_{r})}+ \frac{C}{r-\rho}\|D_{m_0}w\|_{L^2({\cal B}^-_{r})}\\
&&+p^{-p}r^{p(2-\frac{Q+2}{q})+\frac{Q}{2}}\|f\|_{L^{q}({\cal B}^-_{r})}^{p}
\eeno

{\bf Case II:  $f\in L^q$ and $1<p\leq \frac{q(Q+2)}{2(Q+2-q)}$.} Now, we have
\beno
\|I'_{43}\|_{L^{2k}({\cal B}^-_r)}&\leq& C\|\varphi pu^{p-1} f\|_{L^{\frac{2Q+4}{Q+4}}(R^{N+1})}\leq Cp\|f\|_{L^{q}({\cal B}^-_{r})}\|w\varphi\|_{L^2({\cal B}^-_{r})}^{\frac{p-1}{p}}r^{(Q+2)(\frac{Q+4}{2Q+4}-\frac1q-\frac{p-1}{2p})}\\
&\leq& \frac{C}{r-\rho}p^{2\frac{p-1}{p}}\|w\varphi\|_{L^2({\cal B}^-_{r})}+p^{-p}r^{p(Q+2)(\frac{2}{Q+2}-\frac1q+\frac{1}{2p})-1}\|f\|_{L^{q}({\cal B}^-_{r})}^{p}
\eeno
Then
\beno
\|w\varphi\|_{L^{2k}}&\leq& \frac{C}{r-\rho}p^{\frac{2}{1-\tilde{\beta}}}\|w\|_{L^2({\cal B}^-_{r})}+ \frac{C}{r-\rho}\|D_{m_0}w\|_{L^2({\cal B}^-_{r})}\\
&&+p^{-p}r^{p(2-\frac{Q+2}{q})+\frac{Q}{2}}\|f\|_{L^{q}({\cal B}^-_{r})}^{p}.
\eeno

%{\bf Case III: $f\in L^{Q+2}$.} We have
%\beno
%\|w\varphi\|_{L^{2k}}&\leq& \frac{C}{r-\rho}\left(p^{\frac{1}{1-\tilde{\beta}}}\|w\|_{L^2({\cal B}^-_{r})}+ \|D_{m_0}w\|_{L^2({\cal B}^-_{r})}\right)+r^{(Q+2)/2}\|f\|_{L^{q}({\cal B}^-_{r})}^{p}.
%\eeno

Hence, we can complete the proof of  (\ref{eq:Sobolev estimate-w}).   $\diamondsuit$

We also obtain the boundedness property of  nonnegative weak sub-solution of (\ref{eq:general equ}) by the
same argument used by Cinti, Pascucci and Polidoro in \cite{CPP} by using the
Moser's iterative method. We obtain
what follows.

\begin{lemma}[$L^{\infty}$ estimate]\label{lem:infty estimate}
Under the assumptions $({\bf H_1}-{\bf H_3})$,
 let $u$ be a non-negative weak sub-solution of (\ref{eq:general equ}) in $\Omega$.
Let $(x_0,t_0)\in \Omega$ and $\overline{{\cal
B}^-_r(x_0,t_0)}\subset \Omega$ and $p \geq 1$. Then there exists a
positive constant $C$ which depends only on $q,N,\la,Q,$ $\|b'\|_{L^{Q+2}({\cal B}^-_{r}(x_0,t_0))}$, $\|c\|_{L^{q}({\cal B}^-_{r}(x_0,t_0))}$ and $\|f\|_{L^{q}({\cal B}^-_{r}(x_0,t_0))}$ such that,
for $0 < r\leq \frac14$
$$
\sup_{{\cal B}^-_{\fr r 2}(x_0,t_0)} u^p \leq \fr
{C}{r^{Q+2}}\int_{{\cal B}^-_r(x_0,t_0)} u^p+C,
$$
provided that the last integral is finite.
\end{lemma}

{\bf Proof.}
Since $u$ is a non-negative weak sub-solution of (\ref{eq:general equ}), we have
for any nonnegative $\varphi \in C^{\infty}_0(\Omega)$, there holds
\beno
\int_{\Omega} \varphi Yu-(Du)^T AD\varphi \geq\int_{\Omega} \varphi(b'\cdot D_{m_0}u+ cu+f).
\eeno
Without loss of generality, we assume that $\Omega={\cal B}^-_{1}$. Taking $\varphi=\eta(\|x\|)pu^{2p-1}$ and $w=u^p$, where $\eta(x)=1$ if $\|x\|\leq \tilde{\rho}r$ and $\eta(x)=0$ if $\|x\|\geq \bar{\rho}r$ with $\frac12\leq \tilde{\rho}<\bar{\rho}<1$, then we have
\beno
&&\int_{t_0}^{0}\int \left[\frac{1}{2}\eta^2\partial_t (w^2)+ \eta^2(Dw)^T AD w+  2\eta w(Dw)^T AD\eta  \right]dxdt\\
&\leq& -\int_{t_0}^{0}\int\left[\frac12w^2\langle x,BD\rangle(\eta^2)+\eta^2 wb'\cdot D_{m_0}w+ cpw^2\eta^2+pfwu^{p-1}\eta^2\right]dxdt\\
&\leq &-\int_{t_0}^{0}\int\frac12w^2\langle x,BD\rangle(\eta^2)dxdt+I_1+I_2+I_3
\eeno
and $t_0\in (-\bar{\rho}^2r^2,-\tilde{\rho}^2r^2)$ such that
\beno
\int \eta^2w^2(x,t_0)dx\leq \frac{1}{(\bar{\rho}^2-\tilde{\rho}^2)r^2}\int_{-\bar{\rho}^2r^2}^{-\tilde{\rho}^2r^2}\int \eta^2w^2(x,t)dxdt.
\eeno
Moreover, by scaling
\beno
|\langle x,BD\rangle\eta|\leq \frac{C}{(\bar{\rho}-\tilde{\rho})^2r^2},\quad |D_{m_0}\eta|\leq \frac{C}{(\bar{\rho}-\tilde{\rho})r},
\eeno
then
\beno
I_1\leq  \frac12\int_{t_0}^{0}\int\eta^2(Dw)^T AD w+C\|b'\|_{L^{s}({\cal B}^-_{r})}^{\frac2s}\|\eta w\|_{L^{2k}({\cal B}^-_{1})}^2
\eeno
where $\frac1k+\frac2s=1$ and $s=Q+2.$
\beno
I_2\leq p\|c\|_{L^{q}({\cal B}^-_{r})}\|\eta w\|_{L^{2q'}({\cal B}^-_{1})}^2,
\eeno
and
\beno
I_3&\leq& Cp\|f\|_{L^{q}({\cal B}^-_{r})}\|\eta w\|_{L^{2q'}({\cal B}^-_{1})}^{2-\frac{1}{p}}r^{\frac{Q+2}{2pq'}}\\
&\leq&C\|f\|_{L^{q}({\cal B}^-_{r})}(p^4\|\eta w\|_{L^{2q'}({\cal B}^-_{1})}^2+ p^{-2p}r^{\frac{Q+2}{q'}})
\eeno
where $\frac{1}{q}+\frac1{q'}=1$ and obviously $q'<k$.

Let $c_1(r)=\|c\|_{L^{q}({\cal B}^-_{r})}+\|f\|_{L^{q}({\cal B}^-_{r})}$ and $\gamma=1-\frac{Q+2}{2q}\in (0,1)$, then $\frac{Q+2}{q'}=Q+2\gamma$.
Concluding the above estimates, for $t\in (-\frac12r^2,0)$, we get
\ben\label{eq:energy-w}
&&\int \eta^2w^2(x,t)dx+\int_{-r^2}^0\int\eta^2(Dw)^T AD wdxdt\leq \frac{C}{(\bar{\rho}-\tilde{\rho})^2r^2}\int\int \eta^2w^2(x,t)dxdt\nonumber\\
&&+C c_1(r)r^{\frac{Q+2}{q'}} p^4\left(r^{-Q-2} \int_{-r^2}^0\int|\eta w|^{2q'}dxdt\right)^{\frac{1}{q'}}+C c_1(r)r^{\frac{Q+2}{q'}}p^{-2p}
\een

Note that under scaling $u(rx,r^2t)$ still is a solution of (\ref{eq:general equ}) in  ${\cal B}^-_{1}$ with the coefficients $(a_{ij}(rx,r^2t)$, $rb'_i(rx,r^2t)$, $r^2c(rx,r^2t)$ and $r^2f(rx,r^2t)$.
Thus, using the embedding inequality (\ref{eq:Sobolev estimate-w}), for any $\frac12\leq \rho<\tilde{\rho}\leq 1$ we get
\ben\label{eq:Sobolev estimate-w'}
\| w\|_{L^{2k}({\cal B}^-_{\rho r})}&\leq& \frac{C}{\tilde{\rho}-\rho}(p^{\frac{2}{1-\tilde{\beta}}}r^{-1}\|w\|_{L^2({\cal B}^-_{\tilde{\rho} r})}+ \|D_{m_0}w\|_{L^2({\cal B}^-_{\tilde{\rho} r})})\nonumber\\
&&+ Cp^{-p}r^{p(2-\frac{Q+2}{q})+\frac{Q}{2}}\|f\|_{L^q({\cal B}^-_{\tilde{\rho}r})}^p
\een
where $\tilde{\beta}=\frac{Q+2}{q}-1\in (0,1)$ with $q>(Q+2)/2$.

Consequently, it follows from (\ref{eq:energy-w}) and (\ref{eq:Sobolev estimate-w'}) that
\beno
&&\left(r^{-Q-2} \int_{{\cal B}^-_{\rho r}}| w|^{2k}dxdt\right)^{\frac{1}{k}}\leq \frac{C}{(\bar{\rho}-\rho)^4r^{Q+2}}p^{\frac{4}{1-\tilde{\beta}}}\int_{{\cal B}^-_{\bar{\rho} r}} w^2(x,t)dxdt\\
&&+ \frac{C}{(\bar{\rho}-\rho)^2}c_1(r)r^{2\gamma} p^4\left(r^{-Q-2} \int_{{\cal B}^-_{\bar{\rho} r}}|w|^{2q'}dxdt\right)^{\frac{1}{q'}}\\
&&+\frac{C}{(\bar{\rho}-\rho)^2} c_1(r)r^{2\gamma}p^{-2p}+Cr^{2p(2-\frac{Q+2}{q})}p^{-2p}\|f\|_{L^{q}({\cal B}^-_{r})}^{2p}
\eeno
where $\tilde{\beta}=\frac{Q+2}{q}-1\in (0,1)$ and we choose $\tilde{\rho}=\frac{\rho+\bar{\rho}}{2}$.

Since $q'<k$, we let
\beno
p=\left(\frac{k}{q'}\right)^{\nu},\quad \nu=1,2,\cdots,
\eeno
and
\beno
\rho=\rho_\nu=(1+2^{-\nu})\frac12,\quad \bar{\rho}=\rho_{\nu-1}
\eeno
Noting that $(\bar{\rho}-\rho)^4+p^4\leq C^\nu$, following the same arguments as in \cite{Kru64} we have
\beno
\Theta_{\nu}&\doteq&\left(r^{-Q-2}\int_{ {\cal B}^-_{\rho_\nu} } |u|^{2pk}dxdt\right)^{\frac{1}{kp}}\\
&=&\left(r^{-Q-2}\int_{ {\cal B}^-_{\rho_\nu} } |u|^{2k\left(\frac{k}{q'}\right)^{\nu}}dxdt\right)^{\frac{(q')^\nu}{k^{\nu+1}}}\\
&\leq &\left(  C^\nu\Theta_{\nu-1}^p \right)^{\frac1p}+  C^{\frac{\nu}{p}}p^{-2}\\
&\leq & C^{\sum_{\nu\geq 1}\nu \left(\frac{q'}{k}\right)^{\nu} }\Theta_{0}+\sum_{\nu\geq 1}\left[C^{\sum_{\nu\geq 1}\nu \left(\frac{q'}{k}\right)^{\nu} } \left(\frac{q'}{k}\right)^{2\nu} \right]
\eeno
which yields the required result.

In fact, the above analysis shows that $L^\infty$ norm of $u^p$ can be bounded by $L^2$ norm of $u^p$, and one can complete the proof by H\"{o}lder inequality.  $\diamondsuit$

Now we apply Lemma \ref{lem:weak poincare} to the function
$$
w= G(\frac{u}{h}+h^{\fr 18}).
$$
If $u$ is a weak solution of (\ref{eq:general equ}), obviously $w$ is an almost weak
sub-solution as in (\ref{eq:weak subsolution w}). We estimate the value of $I_0$ given by (\ref{eq:I0}) and
(\ref{eq:I1 C2}) in Lemma \ref{lem:weak poincare}.

\begin{lemma}[Lemma 3.4 in \cite{Wang-Zhang2009}]\label{lem:I0 estimate}
Under the assumptions of Lemma \ref{lem:weak poincare}, there exist constants
$\lambda_0$ and $r_0$ with $r_0<\theta$, where $\la_0$ only depends on
constants $\a$,
 $\tilde{\beta}$, $\lambda$, $B$, $N$, and $\varphi$, $0<\lambda_0<1$, such
 that for $r<r_0$ there holds
$$
|I_0|\leq \lambda_0 \ln(h^{-\fr 1 8}).
$$
\end{lemma}

\begin{lemma}\label{lem:lower bound of u}
Suppose that $u(x,t)\geq 0$ be a solution of equation (\ref{eq:general equ}) in
${\cal B}^-_r$ centered at $(0,0)$ and
$$
{\rm meas}\{(x,t)\in {\cal B}^-_r, \quad u \geq 1\} \geq \fr 1 2 {\rm meas} ({\cal
B}^-_r).
$$
Then there exist constant $\theta$ and $h_0$, $0<\theta, h_0<1$
which only depend on $B$, $\la$, $\la_0$ and $N$ such that
$$
u(x,t) \geq h_0\quad \hbox{in}\quad {\cal B}^-_{\theta r}.
$$
\end{lemma}
{\it Proof:} We consider $w=G(\frac{u}{h}+h^{\fr 18})$ for
$0<h<\frac14$, to be decided. {Take $r=\theta h^{\frac{q(Q+2)}{2q-Q-2}},$}  by applying Lemma \ref{lem:weak poincare} to $w$, and we have $$
-\!\!\!\!\!\!\int_{{\cal B}^-_{\theta r}}( w-I_0)_+^2 \leq
C(B,\la)\fr{\theta r^2}{|{\cal B}^-_{\theta r}|} \int_{{\cal
B}^-_r }|D_{m_0}w|^2+C(\theta,B,\la,\|c\|_{L^q},\|u\|_{L^2},\|f\|_{L^{q}})h.
$$
Let $\tilde{u}={\fr u h}$, then $\tilde{u}$ satisfies the conditions
of Lemma \ref{lem:level set}. We can get similar estimates as (\ref{eq:energy estimate})-(\ref{eq:choose of alpha beta}), hence we have
\beno
&& C(B,\la)\fr{\theta r^2}{|{\cal B}^-_{\theta r}|}
\int_{{\cal B}^-_{ r }}|D_{m_0}w|^2\\
& &\leq C(B,\la)\fr{\theta r^2}{|{\cal B}^-_{\theta
r}|}[C(B,\la)(1-{\beta})^{-2}{\beta}^{-Q} +\fr 45\ln(h^{-\fr 18})]
{\rm meas}(K_{\beta r}\times S_{\beta r})\\
& &\leq C(\theta,B,\la) \ln(h^{-\fr 18}),
 \eeno
where $\theta$ has been chosen. By Lemma \ref{lem:infty estimate}, there exists a
constant, still denoted by $\theta$, such that for $z \in {\cal
B}^-_{\theta r}$,
\ben\label{eq: estimate of w-I0}
w-I_0\leq C(B,\la) (\ln(h^{-\fr 18}))^{\fr 12} .
\een
Therefore we may { choose $h_0$ small enough}, so that
$$
C (\ln (\fr {1}{h_0^{\fr 18}}))^{\fr 12}\leq  \ln (\fr
{1}{2h_0^{\fr 18}})-\lambda_0\ln (\fr {1}{h_0^{\fr 18}}).
$$
Then Lemma \ref{lem:I0 estimate} and (\ref{eq: estimate of w-I0}) implies
$$
\max_{{\cal B}^-_{\theta r}}\fr{h_0}{u+h_0^{\fr 98}}\leq \fr
{1}{2h_0^{\fr 18}},
$$
which implies $\min_{{\cal B}^-_{\theta r}}u\geq h_0^{\fr 98}$, then
we finished the proof of this lemma.   $\diamondsuit$

\setcounter{equation}{0}
\section{Proof of our main results}

{\bf Proof of Theorem \ref{thm:main}.} This is similar to that in \cite{Wang-Zhang2009}. We may assume that $M=\max_{{\cal
B}^-_{r}}(+u)=\max_{{\cal B}^-_{r}}(-u)$, otherwise we replace $u$
by $u-C$, since $u$ is bounded locally. Then either $1+\fr u M$ or
$1-\fr u M$ satisfies the assumption of Lemma \ref{lem:lower bound of u}, and we suppose
$1+\fr u M$ does it, thus Lemma \ref{lem:lower bound of u} implies existing $h_0>0$ such
than $\inf_{{\cal B}^-_{\theta r}}(1+\fr u M)\geq h_0,\, $ i.e.
$u\geq M(h_0-1)$,\,then
$$
Osc_{{\cal B}^-_{\theta r}}u\leq M-M(h_0-1)\leq
(1-\fr{h_0}{2})Osc_{{\cal B}^-_{r}}u,
$$
which implies the $C^{\a}$ regularity of $u$ near point $(0,0)$ by
the standard iteration arguments. By the left invariant
translation group action, we know that $u$ is $C^{\a}$ in the
interior.   $\diamondsuit$

\setcounter{equation}{0}
\section{Appendix: G-function}

Next, we introduce some properties of G-function, which was mentioned in \cite{Kru64} (see also \cite{Gu}). Here, we give a detailed description for completeness.

\begin{lemma}[G-function]\label{lem:G function}
There exists a function $G(t):(0,+\infty)\rightarrow \mathbb{R}$ such that\\
%\beno
%&&i)\quad G''(t)\geq G'(t)^2,\quad t>0;\\
%&&ii)\quad G(u)=0,\quad t\geq 1;\\
%&&iii)\quad G'(u)~-\ln t,\quad t\rightarrow 0^+;\\
%&&iv)\quad G(t)\leq \frac{1}{t},\quad t>0.
%\eeno
\begin{align*}\,\, \left\{
\begin{aligned}
&i)\quad G''(t)\geq G'(t)^2,\quad t>0;\\
&ii)\quad G(t)=0,\quad t\geq 1;\\
&iii)\quad G(t)=-\ln t,\quad 0<t\leq \frac14;\\
&iv)\quad  G'(t)\leq 0,\quad t>0.
\end{aligned}
\right. \end{align*}
\end{lemma}

{\it Proof:}
Let $h_0(t)$ be a simple function as follows:
\begin{align*}\,\, h_0(t)=\left\{
\begin{aligned}
&-1,\quad t\leq 1\\
&0,\quad t> 1.
\end{aligned}
\right. \end{align*}
By standard mollifying technique, one can obtain a smooth function $h(t)\in C^{\infty}(\R)$
\begin{align*}\,\, \left\{
\begin{aligned}
&i)\quad h(t)=h_0(t),\quad t\leq \frac12 ~{\rm or}~t>2;\\
&ii)\quad h'(t)\geq 0,\quad t\geq 0;\\
&iii)\quad h(t)\leq 0,\quad t\geq 0;\\
&iv)\quad\int_0^2h(t)dt=-1.
\end{aligned}
\right. \end{align*}
Again, we let $f(t)=\int_0^th(t)dt$, then
\begin{align*}\,\, \left\{
\begin{aligned}
&i)\quad f'(t)=h(t)\leq 0,\quad t\geq 0;\\
&ii)\quad -t\leq f(t)<0,\quad t> 0;\\
&iii)\quad f(0)=0,\quad f(t)=-1,\quad t\geq 2.
\end{aligned}
\right. \end{align*}
Next, write $g(t)=-\ln(-f(t))$, then we have
\beno
g'(t)=-\frac{f'(t)}{f(t)}=-\frac{h(t)}{f(t)}\leq 0,
\eeno
and
\beno
g''(t)=\frac{f'(t)^2-f(t)h'(t)}{f(t)^2}\geq \frac{f'(t)^2}{f(t)^2}=g'(t)^2.
\eeno
Moreover, for $t\leq \frac12$
we have
\beno
g(t)=-\ln(-\int_0^th(t)dt)=-\ln t,\quad 0\leq t\leq \frac12,
\eeno
%\beno
%\lim_{t\rightarrow 0^+}\frac{g(t)}{-\ln t}=\lim_{t\rightarrow 0^+}\frac{g'(t)}{-t^{-1}}=\lim_{t\rightarrow 0^+}\frac{th(t)}{f(t)}=1,
%\eeno
and
\beno
g(t)=-\ln(-\int_0^th(t)dt)=0,\quad t\geq 2.
\eeno
Hence, the proof of (i)-(iii) is complete by choosing $G(t)=g(2t)$.

%Finally, we come to prove (iv). Since the function  $\tilde{g}(t)=g(\mu t+\nu)$ for any $\mu,\nu>0$ satisfies
%\beno
%\tilde{g}'(t)\leq 0,\quad \tilde{g}''(t)\geq 0,\quad for~t\geq 0,
%\eeno
%which implies both $|\tilde{g}'(t)|=-\tilde{g}'(t)$ and $\tilde{g}(t)$ attain its maximum at $t=0$ when $t\geq 0.$ Then
%\beno
%|\tilde{g}'(t)|\leq -\tilde{g}'(0)=-\mu g'(\nu).
%\eeno
%Note that $h=-1$ for $t\leq \frac12$, and we have $-f(t)+th(t)=0$ when $0\leq t\leq \frac12.$ Then for $0<t\leq \frac12$, we get
%\beno
%[-t g'(t)]'&=&-g'(t)-tg''(t)\\
%&=&\frac{h(t)}{f(t)}-t\frac{f'(t)^2-f(t)h'(t)}{f(t)^2}\\
%&\leq &\frac{-h(t)(-f(t)+th(t))}{f(t)^2}=0,
%\eeno
%and
%\beno
%|-t g'(t)|\leq \lim_{t\rightarrow 0^+} |-t g'(t)|=1,
%\eeno
%which yields that
%\beno
%|g'(t)|\leq \frac{1}{t},\quad 0<t\leq \frac12.
%\eeno
Thus, the proof is complete. $\diamondsuit$

\bigskip

\noindent {\bf Acknowledgments.} W. Wang was supported by NSFC under grant 11671067,
 "the Fundamental Research Funds for the Central Universities" and China Scholarship Council. L. Zhang was partially
supported by NSFC under grant 11471320 and 11631008.

\end{document}